\documentclass[11pt]{amsart}
\usepackage{geometry}                
\usepackage{graphicx}
\usepackage{amssymb}
\usepackage{epstopdf}
\usepackage{latexsym}
\usepackage{amssymb}
\usepackage{amsthm}
\usepackage{euscript}

\def\Hu{\H}

\def\Scd{\mathcal{S}(\D)}
\def\Scb{\mathcal{S}(\DT)}

\def\Scn{\mathcal{S}(\D^n)}

\def\={\ = \ }

\def\G{\Gamma}
\def\Ga{\Gamma}

\def\De{\Delta}
\def\La{\Lambda}

\def\l{\lambda}

\def\li{\lambda_i}

\def\vare{\varepsilon}

\let\i=\infty

\def\={ = }    
\def\ot{\otimes}

\def\TT{{\mathbf T}}

\def\Pp{\{\li\to\oo_i\}_1^N}

\def\C{\mathbb C}

\def\T{\mathbb T}
\def\TT{\mathbb T^2} 
\def\D{\mathbb D}
\def\DT{\mathbb D^2} 

\def\be{\setcounter{equation}{\value{theorem}} \begin{equation}}
\def\ee{\end{equation} \addtocounter{theorem}{1}}
\def\beq{\begin{eqnarray*}}
\def\eeq{\end{eqnarray*}}

\def\bp{{\sc Proof: }}
\def\ep{{}{\hfill $\Box$} \vskip 5pt \par}

\def\bl{\begin{lemma}}
\def\el{\end{lemma}}
\def\bt{\begin{theorem}}
\def\et{\end{theorem}}
\def\bprop{\begin{prop}}
\def\eprop{\end{prop}}
\def\bd{\begin{definition}}
\def\ed{\end{definition}}
\def\br{\begin{remark}}
\def\er{\end{remark}}
\def\bexer{\begin{exercise}}
\def\eexer{\end{exercise}}

\newtheorem{theorem}{Theorem}[section]
\newtheorem{prop}[theorem]{Proposition}
\newtheorem{lemma}[theorem]{Lemma}
\newtheorem{cor}[theorem]{Corollary}

\newtheorem{conjecture}[theorem]{Conjecture}
\newtheorem{definition}[theorem]{Definition}
\newtheorem{question}[theorem]{Question}

\newtheorem{remark}[theorem]{Remark}

\newcommand{\Frac}[2] {\displaystyle{\frac{#1}{#2}} }

\newcommand{\oo}{\omega}

\title{Hilbert function spaces and the Nevanlinna-Pick problem on the polydisc II}
\author{David Scheinker}

\begin{document}
\maketitle

\begin{abstract}
In \cite{ds_hfs}, a geometric procedure for constructing a Nevanlinna-Pick problem on $\D^n$ with a specified set of uniqueness was established.
In this sequel we conjecture a necessary and a sufficient condition for a Nevanlinna-Pick problem on $\DT$ to have a unique solution. 
We use the results of \cite{ds_hfs} and Bezout's theorem to establish three special cases of this conjecture.
\end{abstract}

\section{Overview}

The {\em Schur class} of the n$-$disc, $\Scn$, is the set of analytic functions mapping $\D^n$ to $\overline\D$, i.e. satisfying 
$||F||_\i=\sup_{z\in\D^n}|f(z)|\leq1$. 
The Nevanlinna-Pick problem on $\D^n$ is to determine, given distinct nodes $\l_1,...,\l_N\in\D^n$ and target points $\oo_1,...,\oo_N\in\D$, 
whether there exists a function $F\in\Scn$  that satisfies $F(\li)=\oo_i$ for each $i$.  
We are primarily interested in the following question.

\begin{question} What are necessary and sufficient conditions for a Nevanlinna-Pick problem on $\D^2$ to have a unique solution? 
\end{question}

Various authors have studied the uniqueness properties of the Nevanlinna-Pick problem: 
in \cite{baltre98} Ball and Trent show how to parameterize the set of all solutions associated to a given problem on $\D^2$;
in \cite{agmc_three} Agler and McCarthy classify those 2 and 3 point problems on $\D^2$ that have a unique solution; 
in \cite{Kn07} Knese gives sufficient conditions for a 4 point problem on $\D^2$ to have a unique solution; 
in \cite{GHW} Guo, Huang and Wang give sufficient conditions for a 3 point Pick problem on $\D^3$ to have a unique solution;
in \cite{ds_pl}, the present author gives sufficient conditions for a Nevanlinna-Pick problem on $\D^n$ to have a unique solution;
in \cite{ds_hfs}, the present author gives a geometric procedure for constructing a Nevanlinna-Pick problems on $\D^n$ with a specified set of uniqueness. 

In this work we introduce the notion of a {\em strong Pick set} and a question closely related to Question 1.1. 
To state them we recall that a rational function $f\in\Scn$ is called {\em inner} if $|f|=1$ almost everywhere on $\T^n$ 
and that an irreducible algebraic variety $V\subset \C^n$ is called {\em inner} if it meets $\D^n$ and exits $\D^n$ through the n-torus, 
i.e. $V\cap\D^n\neq\emptyset$ and $V \cap \partial ( \D^n) \subset \T^n$. 

\begin{definition} Given a function $f\in\Scn$ and an inner variety $V\subset \C^n$, we say that $V$ is a 
{\em strong Pick set} for $f$, if each $h\in\Scn$ that equals $f$ on $V\cap\D^n$ equals $f$ on $\D^n$, i.e. if $h|_V=f|_V$, then $h=f$.
\end{definition}

\begin{question} Given a rational inner function $f\in\Scb$ and an inner variety $V$, what are necessary and sufficient conditions for $V$ to be a strong 
Pick set for $f$?
\end{question}

The degree of a rational inner function $f$ on $\D^n$, denoted $\deg(f)$, is the degree of the numerator of $f$ in an irreducible representation. The 
degree of a rational inner function $f$ on $\D^n$ in $z_i$, denoted $\deg_i(f)$, is the degree of such a numerator in $z_i$. 
On $\D$, the answer to Question 1.3 is given by the following corollary of Pick's 1916 results.


\begin{cor} \label{PickCor}
For a polynomial $p$ with zeros given by distinct points in $\D$ and a rational inner function $f$, $V=Z_p$ is a strong Pick set for $f$ if and only if $\deg(f)<\deg(p)$.
\end{cor}

We state Pick's original result and derive Corollary \ref{PickCor} in Section 6. Our conjecture is that Corollary \ref{PickCor} generalizes to $\D^2$. 

\begin{conjecture} 
\label{Conj}
Fix a rational inner function $f$ on $\DT$ and an irreducible inner variety $V=Z_p$.\\
If $\deg_i(f)<\deg_i(p)$ for $i=1,2$, then $V$ is a strong Pick set for $f$.\\
If $\deg_i(f)\geq\deg_i(p)$ for $i=1,2$, then $V$ is not a strong Pick set for $f$.
\end{conjecture}

The main results of this paper are several special cases of Conjecture \ref{Conj}. We do not address the mixed case 
$\deg_1(f)<\deg_1(p)$ and $\deg_2(f)\geq\deg_2(p)$ since there exist such examples where $V$ is and is not a strong Pick set for $f$. We also mention that 
the following partial case of this conjecture was established in \cite{ds_pl}.
\begin{theorem}(Scheinker, \cite{ds_pl}) Fix positive integers $n$ and $N$. There exists a 1-dimensional inner variety $V\subset\C^n$ with the following property. 
$V$ is a strong Pick set for each rational inner function $f$ on $\D^n$ that satisfies $\deg(f)<N$.
\end{theorem}

This paper is organized as follows. 
In the Section 2 we give some background and establish the relationship between Question 1.1 and Question 1.3. 
In the Section 3 we state our main results. 
In Sections 4, 5 and 6 we prove our main results. 

I would like to thank Jim Agler and Hugo Woerdeman for  several very useful conversations about this research. I would also like to thank Kelly Bickel for several 
very useful conversations about a special case of Theorem \ref{thmMon}.


\section{Background}



A Nevanlinna-Pick problem on $\D^n$ is called {\em extremal} if a solution $f$ satisfying $||f||_\infty=1$ exists and 
no solution $h$ satisfying $||h||_\infty<1$ exists. If a problem is not extremal,  then it does not have a unique solution. Indeed, if there exists a solution 
$f$ with $||f||_\infty<1$, then for any polynomial $p$ vanishing on the nodes and any $g\in\Scn$ of sufficiently small norm, $f+ pg$ is a 
solution. 
On $\D$, the condition of being extremal is sufficient for a problem to have a unique solution. The following example shows 
that on $\D^2$, unlike on $\D$,  a Nevanlinna-Pick problem may be extremal and fail to have a unique solution.


\begin{exam}
\label{exam2p}
On $\D^2$, the problem with data $(0,0), (\frac{1}{2},\frac{1}{2})$ and $0, \frac{1}{2}$ is extremal and fails to have a unique solution. 
Let $V=Z_p$ where $p=z-w$.  If $f$ is a solution, then $f_d(z)=f|_V=f(z,z)$ is in $\Scd$, satisfies $f_d(0)=0$, $f_d(1/2)=1/2$ and the classical 
Schwarz lemma implies that $f_d(z)=z$. Thus, all solutions to the problem agree on $V\cap\DT$ and the problem is extremal since 
$||f||_\infty\geq||f_d||_\infty=1$. The solution is not unique since each coordinate function solves.
\end{exam}

Example 2.1 is representative of those extremal Nevanlinna-Pick problems on $\D^2$ that fail to have a unique solution, in a sense made precise by the following three theorems.

\begin{theorem} 
\label{Inn}
(Agler, \cite{ag1}): If a Nevanlinna-Pick problem on $\D^2$ has a solution, then it has a solution that is a rational inner function. \end{theorem}

\begin{theorem} (Agler and McCarthy, \cite{agmc_dv}):
\label{DisVar}
Given an extremal Nevanlinna-Pick problem on $\D^2$, there exists an inner variety $V$ with the property that all solutions agree on 
$V\cap\D^2$.
\end{theorem}

\begin{theorem} (Scheinker, \cite{ds_hfs}):
\label{RegInn}
Given a rational  inner function $f$  and an inner variety $V=Z_p$ there exists a Nevanlinna-Pick problem on $\D^2$ with nodes lying on $V$ such that each 
solution to the problem equals $f$ on $V\cap\DT$.
\end{theorem}

These theorems allow us to demonstrate the relationship between Question 1.1 and Question 1.3. 
Suppose that the problem with data $\l_1,...,\l_N$ and $\oo_1,...,\oo_N$ has a unique solution $f$. Theorem \ref{Inn} 
implies that $f$ is a rational inner function and the proof of Theorem \ref{DisVar} guarantees the existence of an inner variety $V$ 
containing the nodes $\l_1,...,\l_N$. If $g\in\Scb$ equals $f$ on $V$, then $f$ is another solution to the problem and $g=f$ on $\DT$. Thus, $V$ is a 
strong Pick set for $f$. Conversely, suppose that $f$ is a rational inner function and that $V$ is a strong Pick set for $f$. Theorem \ref{RegInn} 
guarantees the existence of with nodes lying on $V\cap\D^2$ with the property that $f$ is a solution and that all 
solutions agree on $V\cap\D^2$. If $g$ is another solution to the problem, then $g$ equals $f$ on $V$ which implies that $g=f$. Thus, 
the problem has a unique solution. 

\section{Statement of main results}

Our first main result allows us to establish several cases of Conjecture \ref{Conj}. It is stated using the inner product and the norm of the Hardy space of $\DT$. 
The Hardy space of $\DT$, denoted $H^2$, is the Hilbert space of analytic functions on $\D^2$ with square summable Taylor coefficients at $(0,0)$ 
and norm and inner product given by the following equivalent formulas 
(we recommend \cite{ampi} for a concise presentation of the pertinent facts about $H^2$). 
For  $f=\sum_0^\infty a_nz^n$ and $g=\sum_0^\infty b_nz^n$ in $H^2$,
\[<f,g>=\sum_0^\infty a_n\overline{b_n}=\int_{\TT} f\overline{g}dm \textrm{ and } ||f||_2^2=\sum_0^\infty |a_n|^2=\int_{\TT} |f|^2dm.\]

\begin{theorem}
\label{HS}
Fix a rational inner function $f$  and an inner variety $V=Z_p$.\\
If for each function $g$ analytic on $\D^2$ such that $pg$ is bounded the inequality 
\[2Re(<f,pg>) < ||pg||_2^2\] holds, then $V$ is a strong Pick set for $f$.
\end{theorem}

Theorem \ref{HS} is somewhat surprising since the norm of the Hardy space on $\DT$ is not equivalent to the infinity norm on $\DT$ in which the Nevanlinna-Pick problem 
is stated. The usefulness of Theorem \ref{HS} is, of course, contingent on the difficulty of showing that the hypothesis holds. 
To demonstrate the applicability of Theorem \ref{HS} we mention that the following result, of independent interest, is an almost immediate corollary.


\begin{theorem}
\label{thmMon}
Fix $f=z_1^{d_1}z_2^{d_2}$, fix an inner variety $V=Z_p$ and write $p$ as a sum of monomials, $p=m_1+...+m_k$. 
If for each $m_k$, $\deg_i(f)<\deg_i(m_k)$ for $i=1$ or $i=2$, then $V$ is a strong Pick set for $f$.
\end{theorem}

To examine the implications of Theorem \ref{thmMon}, one could try proving directly the special case when $f=z_1z_2$ and $p=z_1^2-z_2^2$,
i.e. if $g\in\Scb$ equals $f$ on the zero set of $p$, then $g=z_1z_2$.  
The present author is unaware of a simple proof of this seemingly simple result. 
We prove Theorem \ref{HS} and derive Theorem \ref{thmMon} as a corollary in Section 4.

Our second main result establishes Conjecture \ref{Conj} for {\em regular} rational inner functions, rational inner functions that are continuous on a neighborhood of $\DT$.
\begin{theorem} Fix a regular rational inner $f$ and an irreducible inner variety $V=Z_p$.
\label{thmReg}
If $\deg_i(f)<\deg_i(p)$ for $i=1,2$ and $h$ is a regular rational inner function that satisfies $h|_V=f|_V$, then $h=f$.\\
If $\deg_i(f)\geq\deg_i(p)$ for $i=1,2$, then $V$ is not a strong Pick set for $f$.
\end{theorem}

We prove Theorem \ref{thmReg} in Section 5 using Bezout's Theorem.

Our last main result is a complete classification of extremal minimal Nevanlinna-Pick problems on $\DT$ that have a solution of one variable only. 
One may expect the study of such a problem to reduce trivially to the study of a problem on $\D$. However, this is not the case since, as becomes evident from the proof of 
Theorem \ref{OneVariableForward}, there may exist a non-trivial geometric relationship between the first and second coordinates of the nodes.

\bt
\label{OneVariableForward}
Fix an extremal minimal Nevanlinna-Pick problem that has a solution $f$, a function of $z_1$ only. 
There exists an Blaschke product of one variable $m(\l)$ and inner variety $V=Z_p$ such that $V\cap\DT=\{(\l,m(\l)):\l\in\D\}$ contains the nodes of the problem and all solutions to the problem agree on $V$. Furthermore, one of the following holds.\\
If $\deg_1(f) < \deg_1(p)$, then $V$ is a strong Pick set for $f$ and $f$ is the unique solution. \\
If $\deg_1(f)\geq \deg_1(p)$, then is not a strong Pick set for $f$ and $f$ is not the unique solution.
\et

We prove Theorem \ref{OneVariableForward} in Section 6 by generalizing an argument from \cite{ampi}.

\section{Proof of Theorem \ref{HS}}

In this section we prove Theorem \ref{HS} and derive Theorem \ref{thmMon} as a corollary. 


{\sc Proof of Theorem~\ref{HS}:}
Fix a rational inner function $f$, an inner variety $V=Z_p$ and suppose that for each function $g$ analytic on $\D^2$ such that $f-pg$ is bounded the following inequality holds:
\[2Re(<f,pg>) < ||pg||_2^2.\]

Suppose, towards a contradiction, that there exists $r\in\Scb$ such that $r|_V=f|_V$ and $r\neq f$. 
We first show that there exists a rational inner function $h$ that satisfies $h|_V=f|_V$ and $h\neq f$. 
By Theorem \ref{RegInn} there exists a Nevanlinna-Pick problem with nodes $\l_1,...,\l_N\in V\cap\DT$ and target values $f(\l_1),...,f(\l_{N})$ 
such that all solutions agree on $V\cap\DT$. Since $r\neq f$, there exists a $\l_{N+1}\in \D^2$ such that $f(\l_{N+1})\neq r(\l_{N+1})$. Consider the 
Nevanlinna-Pick problem with nodes $\l_1,...,\l_N, \l_{N+1}$ and target values $r(\l_1),...,r(\l_{N+1})$. The problem is solvable since 
$r$ is a solution and Theorem \ref{Inn} implies that there exists a rational inner solution $h$. But notice, since $h$ is also a solution to the problem 
with data $\l_1,...,\l_N$ and $f(\l_1),...,f(\l_{N})$, $h$ equals $f$ on $V\cap\DT$.

Write $f-h=pg$ where $g$ is analytic on $\D^2$ and notice that $pg$ is bounded on $\D^2$ since it 
is the difference of two bounded functions. 
\begin{eqnarray}
 1&=&||h||_\infty^2\\
&=&||f-pg||_\infty^2\\
&=&\int_{\TT} |f-pg|^2 dm\\
&=&||f-pg||_2^2\\
&=&||f||_2^2-2Re<f,pg>+||pg||_2^2\\
&=&1-2Re<f,pg>+||pg||_2^2
\end{eqnarray}
Thus, $2Re(<f,pg>)=||pg||_2^2$ which contradicts our assumption. 
The equality of 4.2 and 4.3 follows from the fact that $h=f-pg$ is inner, i.e. has modulus equal to 1 almost everywhere on $\TT$.
\ep

We now prove Theorem \ref{thmMon}. Fix a rational inner function $f=z_1^{d_1}z_2^{d_2}$ and an inner variety $V=Z_p$. 
Write $p$ as the sum of monomials $p=m_1+...+m_k$ ordered so that for $j=1,...,l$, $\deg_1(f)<\deg_1(m_j)$ and for $j=l+1,...,k$,  $\deg_2(f)<\deg_2(m_j)$.
 

\begin{eqnarray}
<f,pg> & = & <z_1^{d_1}z_2^{d_2},(m_1+...+m_k)g>\\
&=&\sum_{j=1}^k<z_1^{d_1}z_2^{d_2},m_jg>\\
&=&\sum_{j=1}^k \int_{\TT} z_1^{d_1}z_2^{d_2}\overline{m_jg}dm\\
&=&\sum_{j=1}^l \int_{\TT} z_2^{d_2}\overline{z_1^{-d_1}m_jg}dm+\sum_{j=l+1}^k \int_{\TT} z_1^{d_1}\overline{z_2^{-d_2}m_jg}dm\\
&=&\sum_{j=1}^l<z_2^{d_2},z_1^{-d_1}m_jg>+\sum_{j=l+1}^k<z_1^{d_1},z_2^{-d_2}m_jg>\\
&=&0.
\end{eqnarray}
Thus, $<f,pg>=0$ and the conclusion follows from Theorem \ref{HS}. The equality of 4.9 and 4.10 follows from the fact that on $\TT$ one has 
$z_i^{a}=\overline{z_i^{-a}}$. The equality  $<z_2^{d_2},z_1^{-d_1}m_jg>=0$ for $j=1,...,l$ follows from noticing that the 
Taylor coefficients of the $z_2^{d_2}$ term in the Taylor series of $z_1^{-d_1}m_jg$ is zero and the equality $<z_1^{d_1},z_2^{-d_2}m_jg>=0$ 
follows by an analogous consideration. 






\section{Proof of Theorem \ref{thmReg}}

The second part of Theorem \ref{thmReg} is an immediate consequence of the following result. 
\bt
\label{RegularLocal}(Scheinker \cite{ds_hfs})
Let $f$ be a regular rational inner function and $V=Z_p$ an inner variety. If $\deg_i(f)\geq\deg_i(p)$ for $i=1,2$, then there exists a rational inner 
function $g$ that equals $f$ on $V$ and does not equal $f$ on $\DT$.
\et

To prove the first of part Theorem \ref{thmReg} consider a regular rational inner function $f$ with $\deg_i(f)=d_i$ and an inner variety 
$V=Z_p$ with $\deg_i(p)=n_i$ such that $d_1<n_1$ and $d_2<n_2$. 
Let $g$ be a regular rational inner function with $\deg_i(g)=e_i$ such that $g=f$ on $V$. 
Assume, towards a contradiction, that $g\neq f$ on $\DT$. 

The way we have set things up, Theorem 2.8 of \cite{agmc_dv} implies that the number of zeros of $f$ on $V$ is $d_1n_2+d_2n_1$ 
and the number of zeros of $g$ on $V$ is $e_1n_2+e_2n_1$. 
The assumption that $g=f$ on $V$ implies that $d_1n_2+d_2n_1=e_1n_2+e_2n_1$. 
If we can show that $f$ and $g$ have at most $d_{1}e_{2}+d_{2}e_{1}$ common zeros in $\C^2$, then we will have the following inequality
\[e_1n_2+e_2n_1=| Z_f \cap Z_g \cap Z_p | \leq |Z_f\cap Z_g |= e_{1}d_{2}+e_{2}d_{1},\]
which contradicts the assumption that $d_1<n_1$ and $d_2<n_2$.

Thus, the proof of Theorem \ref{thmReg} will be complete once we establish the following theorem.
\bt
\label{InnerIntersection}
Let $f$ and $g$ be rational inner functions of degree $(d_1,d_2)$ and $(e_1,e_2)$. The number of common zeros of $f$ and $g$ counted with multiplicity 
is less than or equal to $d_{1}e_{2}+d_{2}e_{1}$. That is, $|Z_{f}\cap Z_{g}|\leqslant d_{1}e_{2}+d_{2}e_{1}$.
\et

For the reader's convenience we now recall the definitions and results used to state Bezout's Theorem and prove Theorem \ref{InnerIntersection}. 
Rather than discuss the notion of a general algebraic variety in $\C^n$ given as the intersection of the zero sets of several polynomials, 
we specialize the presentation of \cite{Griff} to emphasize the notion of a plane algebraic curve, an algebraic variety in $\C^2$ given as the zero set of a single polynomial.  
To simplify notation and keep with the notation of \cite{Griff} and we use the variables $(x,y)$ instead of $(z_1,z_2)$.

\begin{definition} (I.8.1 \cite{Griff})
An {\em affine algebraic curve} is a subset of $\C^2$ defined by 
\[ V=\{(x,y)\in\C^2: p(x,y)=0\},\] 
where $p$ is a polynomial. The {\em degree} of $V$ is the degree of $p$. We write $V=Z_p$.
\end{definition}

\begin{definition} (I.8.2 \cite{Griff})
A {\em projective algebraic curve} is a subset of $\mathbb{P}^2\C$  defined by 
\[ V=\{(x,y,z)\in\mathbb{P}^2\C: P(x,y,z)=0\},\] 
where $P$ is a homogeneous polynomial. The {\em degree} of $V$ is the degree of $P$. We write $V=Z_P$.
\end{definition}

We use the natural embedding of $\C^2$ into $\mathbb{P}^2\C$ that identifies the points $(x,y)\in\C^2$ and $(x,y,1)\in\mathbb{P}^2\C$. 
We will abuse notation and  write $(x,y)\in\C^2\subset\mathbb{P}^2\C$. 
This identification allows us to identify an affine algebraic curve $V=Z_p$ of degree $n$ with the projective algebraic curve $V=Z_P$ of degree $n$ as follows. 
Given $p(x,y)$ of degree $n$, let $P(x,y,z) =  z^np(\frac{x}{z},\frac{y}{z})$ and given $P(x,y,z)$, let $p(x,y)=P(x,y,1)$.

\bl (II.5.1 \cite{Griff})
\label{Coord}
Suppose $V$ is a projective algebraic curve with $\l \in V$. There exists a coordinate system such that 
\[\l=(0,0)\in\C^2\subset\mathbb{P}^2C,\]
and such that the affine equation of $V$, given by $p(x,y)=0$, satisfies 
\[p(x,y)=y^k+a_1(x)y^{k-1}+...+a_k(x),\]
where $a_j(x)$ is a polynomial of degree less than or equal to $j$ or $a_j(x)=0$.
\el

\begin{definition} (II.7.3 \cite{Griff})
\label{IntNum}
Suppose the affine algebraic curves $V=Z_p$ and $W=Z_q$ intersect at the point $\l$. 
After a suitable change of coordinates, we may assume that $\l=(0,0)$ and that the conclusion of Lemma \ref{Coord} holds.
If $p$ is locally irreducible in a neighborhood of $(0,0)$, then there exists a local normalization of $V$ at $(0,0)$ given by $g:\D\to\DT$ with 
\[g(t)=(t^k,y_v(t^k)),\] 
and we define the {\em intersection number of $V$ and $W$ at $\l=(0,0)$} as the multiplicity of the zero of the one variable analytic function $h(g(t))$ at $t=0$.
In the general case, suppose that in a neighborhood of $\l=(0,0)$, $p$ factors as $=p_1^{m_1}\cdot ... \cdot p_l^{m_l}$ where each $p_j$ is locally irreducible 
in a neighborhood of $\l=(0,0)$. Let $V_j=Z_{p_j}$ and define the intersection number of $V$ and $W$ at $\l=(0,0)$ as
\[(V\cdot W)_{\l}=\sum_{j=1}^l m_j(V_j\cdot W)_{\l}.\]
\end{definition}


\begin{definition}  (II.7.4 \cite{Griff})
The intersection number of two projective algebraic curves $V$ and $W$ in $\mathbb{P}^2\C$ is
\[(V\cdot W)=\sum_{\l\in V\cap W} (V\cdot W)_{\l}.\]
\end{definition}

\bt{\em (II.7.5 Bezout)}\\
Suppose two projective algebraic curves $V=Z_P$ and $W=Z_Q$ have no common curve components (i.e. the polynomials $P$ and $Q$ have no common factor). 
Then 
\[(V\cdot W)=\deg(V)\cdot\deg(W)=\deg(P)\cdot\deg(Q).\]
\et 

Finally, consider two rational inner functions as the ratios of irreducible polynomials $f_1=\frac{q}{q_d}$, $f_2=\frac{r}{r_d}$ and an inner variety $V=Z_p$ with $p$ irreducible.
Let $Q, S$ and $P$ denote the projective polynomials associated to $q,s$ and $p$. 
Define the {\em number of common zeros of $f_1$ and $f_2$} as the sum of the intersection numbers of the projective curves $Z_Q$ and $Z_S$ at points $\l\in\C^2$, i.e. 
\[|Z_{f_1}\cap Z_{f_2}|=\sum_{\l\in(Z_Q\cap Z_S)\cap\C^2}(Z_Q\cap Z_S)_\l.\]
Define the {\em number of zeros of $f$ on $V\cap\DT$} as the intersection numbers of the projective curves $Z_Q$ and $Z_P$ at points in $\D^2$, i.e. 
\[\deg_V(f)=\sum_{\l\in(Z_Q\cap Z_P)\cap\DT}(Z_Q\cap Z_P)_\l.\]




{\sc Proof of Theorem~\ref{InnerIntersection}:}

Write $f_1$ and $f_2$ as the ratios of two polynomials relatively prime in $\C[z,w]$.
\[ f_1(x,y)=\frac{q(x,y)}{q_d(x,y)} \hspace{15pt} \textrm{ and }\hspace{15pt}  f_2(x,y)=\frac{r(x,y)}{r_d(x,y)}\]
If $\overline{a_n}$ and $\overline{b_m}$ are the, necessarily non-zero, constant terms of $q_d$ and $p_d$, then Rudin's theorem on the structure of rational inner functions \cite{rud69} 
implies that $p$ and $q$ have the form
\[q(x,y) =  a_{n}x^{d_{1}}y^{d_{2}} +  a_{n-1}x^{d_{1}-1}y^{d_{2}}... + a_{0} \textrm{ and } r(x,y) =  b_{m}x^{e_{1}}y^{e_{2}} + b_{m-1}x^{e_{1}-1}y^{e_{2}}... + b_{0} \] 

Letting $n=d_1+d_2$ and $m=e_1+e_2$, the projective polynomials corresponding to $q$ and $r$ have the form
\[Q(x,y,z) =  z^n q(\frac{x}{z},\frac{y}{z}) =a_{n}x^{d_{1}}y^{d_{2}} +  a_{n-1}x^{d_{1}-1}y^{d_{2}}z... + a_{0}z^n \] 
\[R(x,y,z) =  z^m r(\frac{x}{z},\frac{y}{z})=b_{m}x^{e_{1}}y^{e_{2}} + b_{m-1}x^{e_{1}-1}y^{e_{2}}z... + b_{0}z^m\] 

Bezout's theorem tells us that the intersection number is 
\[\deg(Q)\cdot\deg(R)=(d_{1}+d_{2})(e_{1}+e_{2})=d_{1}e_{1}+d_{2}e_{2}+d_{1}e_{2}+d_{2}e_{1}.\]

The intersection number of these polynomials at infinity is at points in $\mathbb{P}^2\C$ of the form $\{x,y,0\}$. At these points the polynomials take the form
\[Q(x,y,0) =  a_{n}x^{d_{1}}y^{d_{2}} \textrm{ and } R(x,y,0) =  b_{m}x^{e_{1}}y^{e_{2}}.\] 

Breaking these up gives
\begin{eqnarray*}
Q(x,y,0) =  a_{n}x^{d_{1}}& \textrm{ and }& R(x,y,0) =  b_{m}x^{e_{1}} \textrm{ at points of the form } \{x,1,0\}. \\
Q(x,y,0) =  a_{n}y^{d_{2}} &\textrm{ and }& R(x,y,0) =  b_{m}y^{e_{2}} \textrm{ at points of the form } \{1,y,0\}.
\end{eqnarray*}
These intersect  at $\{0,1,0\}$ with multiplicity $d_{1}e_{1}$ and at $\{1,0,0\}$ with multiplicity $d_{2}e_{2}$.

Subtracting the $d_{1}e_{1}+d_{2}e_{2}$ intersections at infinity from the intersection number gives $d_{1}e_{2}+d_{2}e_{1}$ 
as an upper bound for the number of intersection points of the form $\{x,y,1\}$, i.e. in $\C^2$. 
Thus, $|Z_{f_1}\cap Z_{f_2}|\leq d_{1}e_{2}+d_{2}e_{1}$.
\ep

\section{Proof of Theorem \ref{OneVariableForward}}
In this section we characterize problems that have a solution of one variable only by generalizing an argument from Chapter 12 of \cite{ampi}. 
In the remainder of this section we use $\{\li\to \oo_i\}_1^N$ to denote the Nevanlinna-Pick problem with data $\l_1,...,\l_N$ and $\oo_1,...,\oo_i$. 
Before proving Theorem \ref{OneVariableForward}, we recall several definitions and results and prove Corollary \ref{PickCor} from the introduction.

\begin{theorem}{(Pick 1916)}
\label{Pick1}
On $\D$, the following are equivalent.\\
\textbf{a.} The problem $\{\l_i\to\oo_i\}_1^N$ is solvable.\\
\textbf{b.} The Pick matrix $P =  \displaystyle{ \left( \frac{ 1- \overline{\omega_i}\omega_j  } {  1-\overline{\l_i}\l_j    } \right)}$ is positive semi-definite.\\
\textbf{c.} The problem $\{\l_i\to\oo_i\}_1^N$ has a rational inner solution $f$ with $\deg(f)=rank(P)$.\\
In this case, the following are equivalent.\\
\textbf{i.} The problem has a unique solution.\\
\textbf{ii.} The problem is extremal.\\
\textbf{iii.} The Pick matrix  $P$ is singular.
\et

{\sc Proof of Corollary~\ref{PickCor}:} Fix a rational inner function $f$ on $\D$, fix $V=Z_p$ where $p$ is a polynomial with distinct zeros $\l_1,...,\l_N\in\D$. 
Consider the problem $\{\li\to f(\l_i)\}_1^N$ and the associated Pick matrix, 
\[P =  \displaystyle{ \left( \frac{ 1- \overline{f(\l_i)}f(\l_j)  } {  1-\overline{\l_i}\l_j    } \right)}\] 
Parts $c$ and $iii$ of Theorem \ref{Pick1} imply that the problem has a unique solution if and only if $\deg(f)=rank(P)<N=\deg(p)$. 
Notice that a function $g\in\Scd$ is a solution if and only if $g|_V=f|_V$. 
If $\deg(f)<\deg(p)$, then the problem has a unique solution and each $g\in\Scd$ that satisfies $g|_V=f|_V$ must equal $f$. 
If $\deg(f)\geq\deg(p)$, then the problem fails to have a unique solution and there exists a $g\in\Scd$ such that $g|_V=f|_V$ and $g\neq f$. \ep


Given a a problem $\Pp$, write $\l_i=(\l_i^1,\l_i^2)$ and let $W, \, \La^1$ and $\La^2$  denote the following $N$-by-$N$ matrices.
\[W \= \left( 1 - \bar w_i w_j \right)_{i,j=1}^N \hspace{20pt}
\La^1 \= \left( 1 - \bar \l^1_i \l^1_j \right)_{i,j=1}^N \hspace{20pt}
\La^2 \= \left( 1 - \bar \l^2_i \l^2_j \right)_{i,j=1}^N 
\]
For a matrix $A$, write $A\geq0$ if $A$ is positive semi-definite and $A>0$ if it is positive definite. 
Let $W\cdot K=(W_{ij}K_{ij})$ denote the Schur entrywise product of two matrices $W$ and $K$. 
A positive definite matrix $K$ is an {\em admissible kernel} if $\La^1\cdot K\geq0$ and $\La^2\cdot K\geq0$, and $K$ is {\em active} if $\det(W\cdot K)=0$. 
Finally, if the problem is extremal and no $N-1$ point subproblem $\{\l_{i_k}\to \oo_{i_k}\}_1^{N-1}$ is extremal, 
then the problem is called {\em minimal}.

\bt{(Agler, \cite{ag1})}
\label{Ag1}
On $\DT$, the following are equivalent.\\
i. The problem $\Pp$ has a solution.\\
ii. For each admissible kernel $K$, $W\cdot K\geq0$.\\
iii. There exists a pair of positive semi-definite matrices $(\Ga,\De)$ such that $W=\Ga\cdot \La^1+ \De\cdot\La^2$.
\et

\bl (Agler, McCarthy \cite{agmc_dv})
If $\Pp$ is an extremal Nevanlinna-Pick problem, then there exists an admissible kernel $K$ that is active. 
Furthermore, if the problem is minimal and $K$ is an active kernel, then rank$(K\cdot W) \, = \, N-1$.
\label{lem1}
\el

We now prove the lemmas we need to establish Theorem \ref{OneVariableForward}.

\bl
\label{GenPick}
Fix a problem $\Pp$  on $\DT$ and let $K$ denote the Szeg\Hu{o}  kernel of the Hardy space on $\D^2$,
\[K_{\li,\l_j}=\Frac{1}{(1-\bar \li z_1)(1-\bar \l_j z_2)}.\]
If $W\cdot K$ is singular, then the problem has a unique solution.
\el

\bp
Theorem 1.6 of \cite{ds_hfs} implies that the generalized problem $\Pp$ in the multiplier algebra of $H^2(\DT)$, $Mult(H^2(\DT))$, has a 
unique solution. However, since a multiplier $M_f$ is in the unit ball of $Mult(H^2(\DT))$ if and only if $f\in\Scb$, this implies that the original problem also has a unique solution. \ep

\bl
\label{nonzero}
Fix a problem $\Pp$ on $\DT$. If there exists a pair of non-zero positive semi-definite matrices $(\Ga,\De)$ 
such that $W=\Ga\cdot \La^1+ \De\cdot\La^2$, then there exists a solution $f$ that is a function of both $z_1$ and $z_2$.
\el

\bp
In Theorem \ref{Ag1}, the proof of $iii$ implies $i$  proceeds by showing that the entry wise equalities of $W=\Ga\cdot \La^1+ \De\cdot\La^2$ 
extend to all of $\DT$ in the following sense. There exists a pair of positive semi-definite functions 
$\Ga,\De$ on $\DT\times\DT$ such that $\Ga(\l_i,\l_j)=\G_{ij}$, $\De(\l_i,\l_j)=\De_{ij}$ and a rational inner function $f\in\Scb$ such that
\begin{equation}
\label{AgReal}
\forall(\l,\zeta)\in\DT\times\DT \hspace{20pt} 1 - \overline{f(\l)} f(\zeta) = ( 1 - \overline{\l^1} \zeta^1 )\Ga(\l,\zeta)+ (1 - \overline{\l^2} \zeta^2)\De(\l,\zeta)
\end{equation}

In \cite{ColWer}, Cole and Wermer show that if $f$ is written as the ratio of relatively prime polynomials $\frac{\tilde r}{r}$ then the following version of the Agler realization holds with 
$A_i$ and $B_i$ polynomials.
\begin{equation}
\label{CWReal}
 \overline{r(\l)} r(\zeta) - \overline{\tilde r(\l)}\tilde r(\zeta) = ( 1 - \overline{\l^1} \zeta^1 )\sum_1^M \overline{A_i(\l)}A_i(\zeta) +  (1 - \overline{\l^2} \zeta^2)\sum_1^M \overline{B_i(\l)}B_i(\zeta)
\end{equation}

Suppose, towards a contradiction, that $f$ does not depend on $z_2$. Then, neither $\tilde p$ nor $p$ depends on $z_2$. 
Differentiating both sides of \ref{CWReal} with respect to $\zeta^2$ gives 
\begin{equation}
\label{CWReal2}
0 = ( 1 - \overline{\l^1} \zeta^1 )\frac{d}{d\zeta_2}\sum_1^M \overline{A_i(\l)}A_i(\zeta) + \frac{d}{d\zeta_2}(1 - \overline{\l^2} \zeta^2)\sum_1^M \overline{B_i(\l)}B_i(\zeta)
\end{equation}
Notice that if $\frac{d}{d\zeta_2}\sum_1^M \overline{A_i(\l)}A_i(\zeta)\neq0$, then one can solve for $1 - \overline{\l^1} \zeta^1$ as a ratio of polynomials that 
depend on $\l^2$ and $\zeta^2$, a contradiction. Thus, $\frac{d}{d\zeta_2}\sum_1^M \overline{A_i(\l)}A_i(\zeta)=0$ and (\ref{CWReal2}) can be written as
\[0 = \frac{d}{d\zeta_2}(1 - \overline{\l^2} \zeta^2)\sum_1^M \overline{B_i(\l)}B_i(\zeta).\]
This implies that $\sum_1^M \overline{B_i(\l)}B_i(\zeta)=0$, which implies that $\De(\l,\zeta)=0$  a contradiction. \ep


The following lemma is a slightly modified version of Lemma 12.11 in \cite{ampi}.
\bl
\label{Rank}
Fix an extremal, minimal problem $\Pp$ on $\DT$. If $K$ is an active kernel and $(\Ga,\De)$ is a pair of positive matrices with 
rank$(\Ga)=N-1$ that satisfy $W=\Ga\cdot \La^1+ \De\cdot\La^2$, then rank$(K\cdot \La ^ 1)=1$.
\el

\bp
If $K$ is an active kernel, then $K \cdot W$ has rank $N-1$ and
annihilates some vector $\vec{\gamma} = ( \gamma_1, ... , \gamma_N)$ and since the problem is minimal, each $\gamma_i\neq0$.

Since $\G$ is positive we let $\vec{u^k} = (u_1^k, ... , u_N^k)^T$ and write $\G$ as the sum of rank one matrices $\displaystyle{\G \= \sum_1^{N-1} \vec{u^k} \otimes \vec{u^k} }$
where non of $\vec{u^k}$ and $\vec{u^l}$ are collinear for $k\neq l$ and $\vec{u^k} \otimes \vec{u^k}$ denotes the matrix 
$\left( \vec{u^k} \otimes \vec{u^k} \right)_{ij} \= u_i^k \bar u_j^k$.
Let the rank of $K \cdot \La^1$ equal $P$ and write $\displaystyle{K \cdot \La^1 = \sum_1^{P}\vec{x^r} \otimes \vec{x^r}}$. 
Notice that $\left( K \cdot \La^1 \cdot \Ga \right) \vec{\gamma} =0$, 
since $0=W\cdot K\gamma = (K \cdot \La^1 \cdot \Ga+K \cdot \La^2 \cdot \De)\gamma$ and both matrices on the right are positive semi-definite.

The equality $\left( K \cdot \La^1 \cdot \Ga \right) \vec{\gamma} =0$ implies that all of $\{(\vec{u_k} \otimes \vec{u_k}) \cdot (\vec{x_r} \otimes \vec{x_r})\}_{k=1,r=1}^{N-1,P}$ annihilate 
$\vec{\gamma}$ for each $1\leqslant k\leqslant M$ and $1\leqslant r \leqslant P$, i.e. $\displaystyle{ 0 = \sum_{j=1}^N \bar u_j^k \bar x_j^r \gamma_j }$.
Therefore each of the vectors $\vec{\gamma}\cdot\vec{x^r}=( \overline{x_1^r}{ \gamma_1}, ... ,\overline{x_N^r}{ \gamma_N})^T$
is orthogonal to each of $\vec{u^k}$. That is, the vectors $\{\vec{\gamma}\cdot\vec{x^r}\}_1^P$ are contained in the subspace of $\C^n$ perpendicular to the $N-1$-dimensional 
subspace of $\C^N$ spanned by $\{\vec{u^k}\}_1^{N-1}$, i.e. a subspace of $\C^N$ of dimension $1$. 
As none of the entries of $\vec{\gamma}$ are $0$, the vectors $\{\vec{x^r}\}_1^P$ must all be collinear and the rank of $K\cdot \La^1$ is 1. \ep

{\sc Proof of Theorem~\ref{OneVariableForward}:}


Let 
\[\displaystyle{A=\frac{\La^2}{\La^1}=\left( \frac{ 1-\l_i^2 \overline{\l_j^2}} {  1-\l_i^1 \overline{\l_j^2}   } \right)}.\]
Notice 
that the $A$ is the Pick matrix corresponding to the problem $\{\l_i^1\to \l_i^2\}_1^N$ on $\D$.

We first show that $A$ is positive semi-definite. Since $f(z_1)$ is the solution to the problem $\Pp$, it is the solution to the one variable problem with data $\{\l_i^1\to \oo_i\}_1^N$. 
By Theorem \ref{Pick1} the matrix 
\[ \Ga_0= \left(\frac{W_{ij}}{\La^1_{ij}}\right)=\left( \frac{ 1-\omega_i \overline{\omega_j}  } {  1-\l_i \overline{\l_j}   } \right)\]
is positive semi-definite.
Furthermore, since the problem is extremal and minimal, Theorem \ref{Pick1} implies that $\deg_1(f)=$rank$(\Ga_0)=N-1$. Let $K$ be an active kernel for the original problem. 
By lemma \ref{Rank} rank$(K \cdot \La^1)=1$. 
Since $(K \cdot \La^1)$ is positive semi-definite with non-zero diagonal entries the fact that it has rank 1 implies that all of its entries are non-zero. 
The matrix 
\[(\frac{1}{K\cdot \La^1})_{ij}=\frac{1}{K_{ij}\cdot \La_{ij}^1}\]
is also positive semi-definite. 
We conclude that $A$ is positive semi-definite by writing $A=\La^2 \cdot \frac{1}{\La^1} \= K\cdot \La^2 \cdot \frac{1}{K\cdot \La^1}$ and noticing that the right-hand 
side is a Schur product of two positive matrices.

To construct $V$, let $n_1=$rank$(A)$. By the one dimensional Pick theorem there exists a Blaschke product $m(\l)$ of degree $n_1$ such that $m(\li^1)=\li^2$. Write $m$ as the ratio of two 
irreducible polynomials $m(\l)=\frac{q(\l)}{r(\l)}$, let $p(z_1,z_2)=z_2r(z_1)-q(z_1)$ and notice that $V\cap\D^2=Z_p\cap\DT=\{ (\l, m(\l)) \, :\, \l \in \D \}$ contains all of the nodes 
of the original problem. Furthermore, since $V$ is inner and the restriction of $f$ to $V$ has less than $N$ zeros, each solution of the original problem $\Pp$ equals $f$ on $V$ by 
Theorem 1.7 of \cite{ds_hfs}. We now examine two cases:

\textbf{Case i. $\deg_1(f)<\deg_1(p)$}. To show that $V$ is a strong Pick set for $f$, fix a $g\in\Scb$ that equals $f$ on $V$. 
There exists a point $w\in \D$ such that $V\cap\DT$ contains $n_1$ distinct points of the form $(l_j,w)$. 
Consider the problem $\{(l_j,w)\to f(l_j)\}_1^{n_1}$ on $\DT$ 
and consider the matrix $W\cdot K$ associated to this problem with $K$ the Szeg\Hu{o}  kernel of $\DT$,
\[ W\cdot K=\left(\Frac{1-\overline{f(l_i)}f(l_j)}{(1-\bar z_i z_j)(1-\bar w w)}\right)=
\left(\Frac{1-\overline{f(l_i)}f(l_j)}{(1-\bar z_i z_j)(1-|w|^2)}\right)=
\Frac{1}{(1-|w|^2)}\left(\Frac{1-\overline{f(l_i)}f(l_j)}{(1-\bar z_i z_j)}\right).\]
The right most matrix in the above equality has rank equal to $\deg(f)<n_1$ by Theorem \ref{Pick1}. 
Thus, $W\cdot K$ is singular, Lemma \ref{GenPick} implies that $f$ is the unique solution to the problem 
$\{(l_j,w)\to f(l_j)\}_1^{n_1}$ and since $g$ is another solution, $g=f$.

\textbf{Case ii. $\deg_1(f)\geq \deg_1(p)$}. To show that $V$ is not a strong Pick set for $f$,   it will suffice to construct a solution $g$ to the original problem that does not equal $f$, 
since in the first part of the proof we showed all solutions agree on $V$. To construct such a $g$ we modify the argument in Chapter 12 of \cite{agmc_dv} and show that there exists a pair of positive semi-definite matrices $(\Ga,\De)$ with $\De$ non-zero such that 
\begin{equation}
\label{realize}
W=\Ga\cdot\La^1+\De\cdot\La^2.
\end{equation}
By Lemma \ref{nonzero}, the existence of such matrices implies that the original problem has a solution that depends on $z_2$.

A pair positive semi-definite matrices $(\Ga,\De)$ satisfies \ref{realize} if and only if 
\begin{equation}
\label{PerPair}
A \cdot \De=\frac{\La^2}{\La^1}\cdot\De \leq \Ga_0
\end{equation}
in which case $\Ga = \Ga_0 - A \cdot \De$.
Write $A$ as the sum of $M$ rank one matrices $\displaystyle{ A \= \sum_{i=1}^M [x^i \ot x^i]}$ and notice that the rank one matrix $\De = \vare [ v \ot v]$
will satisfy (\ref{PerPair}) for some $\vare > 0$ if and
only if for each $r$ the vector 
$ v \cdot x^r \ := \ ( v_1 x_1^r, ... , v_N x_N^r)^T$ lies in the
range of $\Ga_0$. Since the rank of $\Ga_0$ is $N-1$, is suffices to fix any non-zero vector $u$ perpendicular to the range of $\Ga_0$ and 
find $v$ so that for each $r$ the vector $v \cdot x^r$ is perpendicular to $u$. These two constraints translate into the following system of $M$ linear equations 
$$\sum_{i=1}^N v_i x_i^r \bar u_i \= 0 \textrm{ for } r=1,...,M.$$
Since $M<N$, there is a non-zero $v$ in $\C^N$ satisfying the above constraints and hence there exists a rank one $\De$ satisfying $A \cdot \De \ \leq \ \Ga_0$.  \ep


\end{document}